\documentclass[11pt,fleqn]{article} %,draft

\usepackage{amsmath,amssymb,amsthm,amsfonts} % ,enumerate ,backref ,refcheck
\usepackage[mathscr]{eucal}

\newcommand{\const}{{\rm const}}

\newtheorem{theorem}{Theorem}%[section]

{\theoremstyle{definition}

}

\begin{document}

\par\noindent {\LARGE\bf
Group Classification (1+2)-dimensional Linear\\
 Equation of Asian Options Pricing
\par}

\vspace{5.5mm}\par\noindent{\large
\large Stanislav V. Spichak$^{\dag}$, Valeriy I. Stogniy$^{\ddag}$ and Inna M. Kopas$^{\ddag}$
}

\vspace{5.5mm}\par\noindent{\it\small
$^\dag$Institute of Mathematics of NAS of Ukraine, Kyiv, Ukraine
\par}

\vspace{2mm}\par\noindent{\it\small	
$^\ddag$National Technical University of Ukraine ``Igor Sikorsky Kyiv Polytechnic Institute'', Kyiv, Ukraine
\par}

\vspace{4mm}\par\noindent
E-mails:
spichak@imath.kiev.ua,
 stogniyvaleriy@gmail.com,
 innakopas5@gmail.com

\vspace{8mm}\par\noindent\hspace*{10mm}\parbox{140mm}{\small\looseness=-1
We consider a class of (1+2)-dimensional linear partial differential equations of Asian options pricing. Special cases have been used to models of financial mathematics. We carry out group classification of a class equations. In particular, the maximum dimension of Lie invariance algebra within the above class is eight-dimensional. It is shown that an equation with such an algebra can be transformed into the linear Kolmogorov equation with the help of the point transformations of variables. Using the operators of invariance algebra symmetry reduction is carried out and invariant exact solutions are constructed for some equations.
}\par\vspace{4mm}

\section{Introduction}\label{sec:Introduction}

Investigation of Asian options leads to interesting mathematical models which can be often formulated in term of linear partial differential equations. Such models can be found in~\cite{spich-stog-1,spich-stog-2}. All these partial differential equations can be grouped into a general class of the form
\begin{gather}\label{Sp-St-1}
\frac{\partial V}{\partial
\tau}+\frac{1}{2}\sigma^2S^2\frac{\partial^2 V}{\partial
S^2}+rS\frac{\partial V}{\partial S}+f(S)\frac{\partial V}{\partial
A}-rV=0 ,
\end{gather}
where $\tau\in [0;T]$, $T$ is the term of the contract; $V=V(\tau,S,A)$ is the function of the option value, $S$ is the value of the underlying asset; $A$ is the average value of all available prices $S$ of the underlying assets by the time $\tau$; $r$ and $\sigma$ are the constants describing the risk-free interest rate and stock volatility respectively, $f(S)\neq \text{const}$ is arbitrary smooth function of the variable $S$. If $f(S)=\text{const}$, then the corresponding equation from the class~\eqref{Sp-St-1} can be reduced into a linear partial differential equation in two independent variables $(f(S)=0)$ by corresponding Galilean transformation.

The class equations~\eqref{Sp-St-1} contain a few well-known equations of Asian options. Thus, in~\cite{spich-stog-1} were considered the equations with $f(S)=S$ and $f(S)=\ln S$, in~\cite{spich-stog-2} were considered the equation with $\displaystyle f(S)=\frac{1}{S}$.

Wide application of equations~\eqref{Sp-St-1} in the problems of financial mathematics causes an undeniable interest in obtaining its exact solutions. One of the most effective methods for constructing exact solutions is the methods of group analysis \cite{spich-stog-3}--\cite{spich-stog-5}. The theoretical- group methods allow you to integrate differential equations, which have the non-trivial invariance group. That is why the urgent task is a complete group classification of differential equations with arbitrary function that allows selecting the equations with broad symmetry properties from the given class of equations.

We first simplify equation~\eqref{Sp-St-1} by the following change of variables (see~\cite{spich-stog-1}):
\begin{gather}\label{Sp-St-2}
\displaystyle u(t,x,y)=x^m e^{qt}V\left(T-\frac{2t}{\sigma^2},x,\frac{2y}{\sigma^2}\right),\ m=\frac{r}{\sigma^2},\ q=m^2+m .
\end{gather}
Then the differential equation is transformed into
\begin{gather}\label{Sp-St-3}
\displaystyle \frac{\partial u}{\partial t}=x^2\frac{\partial^2 u}{\partial x^2}+f(x)\frac{\partial u}{\partial y} .
\end{gather}

Transformation~\eqref{Sp-St-2} is a point transformation and, therefore, equations \eqref{Sp-St-1} and \eqref{Sp-St-3} equivalent.
It is more convenient to find symmetries for equations~\eqref{Sp-St-3}.

Our goal is to carry out of group classification of class~\eqref{Sp-St-3}. An investigation of symmetry properties of the
equation~\eqref{Sp-St-1} with $f(S)=S$ has been considered in~\cite{spich-stog-6,spich-stog-7} and the equation~\eqref{Sp-St-3} with
$f(x)=-x$ has been considered in~\cite{spich-stog-8}.

\section{Algebra $A^{ker}$ and Group of Equivalence Transformations}\label{sec:Algebra-Ker}

Group classification of equations~\eqref{Sp-St-3} will be carried out using the classical Lie-Ovsyannikov method \cite{spich-stog-3,spich-stog-5}.

Consider a one-parameter Lie group of local transformation in space of variables $(t,x,y,u)$ with an infinitesimal operator of form
\begin{gather}\label{Sp-St-4}
\displaystyle X=\xi^0(t,x,y,u)\partial_t+\xi^1(t,x,y,u)\partial_x+\xi^2(t,x,y,u)\partial_y+\eta(t,x,y,u)\partial_u ,
\end{gather}
which keeps equation~\eqref{Sp-St-3} invariant. The Lie criterion of infinitesimal invariance is
\begin{gather}\nonumber
\tilde X(u_t-x^2u_{xx}-f(x)u_y)\biggl|_{u_t=x^2u_{xx}+f(x)u_y}=0 ,
\end{gather}
where $\displaystyle u_t=\frac{\partial u}{\partial t}$, $\displaystyle u_y=\frac{\partial u}{\partial y}$,
$\displaystyle u_{xx}=\frac{\partial^2 u}{\partial x^2}$, $\tilde X$ is second prolongation of operator~\eqref{Sp-St-4}. This criterion yields
the following determining equations for functions $\xi^0$, $\xi^1$, $\xi^2$, $\eta$, and for arbitrary element $f(x)$:
\begin{gather}\label{Sp-St-5}
 \xi^0_x=\xi^0_u=\xi^1_u=\xi^2_x=\xi^2_u=\eta_{uu}=0 ,
\end{gather}
\begin{gather}\label{Sp-St-6}
 x\xi^0_t-2x\xi^1_x+2\xi^1-xf(x)\xi^0_y=0 ,
\end{gather}
\begin{gather}\label{Sp-St-7}
 \xi^1_t-x^2\xi^1_{xx}-f(x)\xi^1_y+2x^2\eta_{ux}=0 ,
\end{gather}
\begin{gather}\label{Sp-St-8}
 f'(x)\xi^1+f(x)(\xi^0_t-\xi^2_y)-(f(x))^2\xi^0_y+\xi^2_t=0 ,
\end{gather}
\begin{gather}\label{Sp-St-9}
 \eta_t-x^2\eta_{xx}-f(x)\eta_y=0 .
\end{gather}
Equations~\eqref{Sp-St-5} do not contain arbitrary element. Integration of them yields
\begin{gather}\label{Sp-St-10}
\xi^0=\xi^0(t,y),\ \xi^1=\xi^1(t,x,y),\ \xi^2=\xi^2(t,y),\ \eta=\alpha(t,x,y)u+\beta(t,x,y),
\end{gather}
where $\alpha(t,x,y)$, $\beta(t,x,y)$ are arbitrary smooth functions. Equations \eqref{Sp-St-6}--\eqref{Sp-St-9} containing arbitrary element explicitly are called classifying equations. The group classification of \eqref{Sp-St-3} reduces to solving classifying conditions with respect to the coefficients of the operator \eqref{Sp-St-4} and arbitrary element simultaneously.

At first, we are finding the kernel algebra $A^{ker}$.

Splitting system \eqref{Sp-St-6}--\eqref{Sp-St-9} with respect to the arbitrary element and their non vanishing derivatives gives the equations
\begin{gather}\nonumber
\xi^0_t=\xi^0_y=\xi^1=\xi^2_t=\xi^2_y=\alpha_t=\alpha_x=\alpha_y=\beta_y=0 ; \ \ \beta_t-x^2\beta_{xx}=0
\end{gather}
for the coefficients of operators \eqref{Sp-St-4} from the algebra $A^{ker}$ of equation \eqref{Sp-St-3}. Integration of them yields
\begin{gather}\label{Sp-St-11}
\xi^0=C_1,\ \xi^1=0,\ \xi^2=C_2,\ \eta=C_3u+\beta(t,x),
\end{gather}
where $C_1$, $C_2$, $C_3$ are arbitrary constants, $\beta(t,x)$ is an arbitrary solution of equation
$\beta_t=x^2\beta_{xx}$. As a result, the following theorem is true.

\begin{theorem}\label{thm:Number1}
For arbitrary function $f(x)$ the Lie symmetry algebra of equation \eqref{Sp-St-1} is
\begin{gather}\label{Sp-St-12}
A^{ker}=<\partial_t,\partial_y,u\partial_u,\beta(t,x)\partial_u> ,
\end{gather}
where $\beta(t,x)$ is an arbitrary solutions of equation $\beta_t=x^2\beta_{xx}$.
\end{theorem}

The problem of group classification consists in finding of all possible inequivalent cases when the solution of system
\eqref{Sp-St-5}--–\eqref{Sp-St-9} leads to Lie algebra $A^{max}$ that satisfies the condition $A^{max}\supset A^{ker}$.

The next step of algorithm of group classification is finding set of equivalence transformations of class \eqref{Sp-St-3}
which form the equivalence group $G^{equiv}$. In order to construct the equivalence group $G^{equiv}$ we are using direct
method proposed in~\cite{spich-stog-5}. Thus, we obtain the following statement.

\begin{theorem}\label{thm:Number2}
The equivalence group $G^{equiv}$ of class \eqref{Sp-St-3} is formed by the transformations
\begin{gather}
 \nonumber
\bar{t}=\epsilon^2_1t+\epsilon_2 ;\quad \bar{x}=\epsilon_3x^{\epsilon_1} ; \quad \bar{y}=\epsilon_4t+\epsilon_5y+\epsilon_6 ; \\
\label{Sp-St-13} \displaystyle \bar{u}=\epsilon_7e^{(1-\epsilon^2_1)t/4}x^{(\epsilon_1-1)/2}u+\varphi(t,x) ; \quad \bar{f}=\frac{1}{\epsilon^2_1}(\epsilon_5f-\epsilon_4) ,
\end{gather}
where $\epsilon_i$, $i\in\{1,\ldots,7\}$ are arbitrary constants, $\epsilon_1\epsilon_3\epsilon_5\epsilon_7\neq 0$,
$\varphi(t,x)$ is an arbitrary solutions of equation
\begin{gather}\nonumber
\varphi_t=x^2\varphi_{xx}+(1-\epsilon_1)x\varphi_x .
\end{gather}
\end{theorem}

Note that these transformations \eqref{Sp-St-13} will be further applied to unite equations \eqref{Sp-St-3} with the same Lie symmetry.

\section{Classification of Lie symmetries}\label{sec:Lie-symmetries}

The next step of algorithm of group classification is finding of a complete set of inequivalent equations (3) with respect
to the transformations from $G^{equiv}$ which are invariant under Lie algebra $A^{max}$ that satisfies the condition
$A^{max}\supset A^{ker}$ ($\dim A^{max}>\dim A^{ker}$). Hereinafter, the notation $A^{max}$ implies the fulfillment of this condition.

To construct all possible forms of the functions $f(x)$  for which the corresponding set of solutions of system
\eqref{Sp-St-5}--\eqref{Sp-St-9} is wider than the “trivial” solution \eqref{Sp-St-11}.

In article \cite{spich-stog-8}, the class of equations containing equations \eqref{Sp-St-3} was considered, and conditions for variable transformations that convert one equation from this class into another equation of the same class were obtained. Using these transformations,
the condition $\xi_y^0=0$ can be obtained for equations of class \eqref{Sp-St-3}. Thus, we will analyse equations
\eqref{Sp-St-5}--\eqref{Sp-St-9} under the condition $\xi_y^0=0$.

Substituting $\xi_y^0=0$ in equations \eqref{Sp-St-6}, we obtain that
\begin{gather}\label{Sp-St-14}
\displaystyle \xi^1=\left(\frac{1}{2}\xi_t^0(\ln x-1)+P(t,y)\right)x,
\end{gather}
where $P(t,y)$ is arbitrary smooth function.

Consider the simplest case when $\xi_t^0=0$ and $P(t,y)=0$. Then the solution \eqref{Sp-St-5}--\eqref{Sp-St-9}
leads to $f(x)=\const$ or leads to Lie algebra $A^{ker}$ for arbitrary function $f(x)$.

Thus, $\xi_t^0\neq 0$ or $P(t,y)\neq 0$.

Given this, we substitute \eqref{Sp-St-14} in \eqref{Sp-St-8} and obtain the differential equation for the function $f(x)$
\begin{gather}\label{Sp-St-15}
\displaystyle f'+\frac{\xi^0_t-\xi^2_y}{\left(\frac{1}{2}\xi^0_t(\ln x-1)+P(t,y)\right)x}f=
-\frac{\xi^2_t}{\left(\frac{1}{2}\xi^0_t(\ln x-1)+P(t,y)\right)x} .
\end{gather}

Since the function $f(x)$  depends only on the variable $x$, the following ordinary differential equation
for the function is obtained from equation \eqref{Sp-St-15}
\begin{gather}\label{Sp-St-16}
\displaystyle f'+\frac{a_1}{(a_2\ln x+a_3)x}f=\frac{a_4}{(a_2\ln x+a_3)x} ,
\end{gather}
where $a_i$, $i\in\{1,\ldots,4\}$ are arbitrary constants.

Solutions of equation \eqref{Sp-St-16} will be of the following forms (solution $f(x)=\const$ are not considered here):

1)\ if $a_1\neq 0$, $a_2=0$, $a_3\neq 0$, then $\displaystyle f=Cx^{-a_1/a_3}+\frac{a_4}{a_1}$;

2)\ if $a_1\neq 0$, $a_2\neq 0$, then $\displaystyle f=C(a_2\ln x+a_3)^{-a_1/a_2}+\frac{a_4}{a_1}$;

3)\ if $a_1=a_2=0$, $a_3\neq 0$, $a_4\neq 0$, then $\displaystyle f=\frac{a_4}{a_3}\ln x+C$;

4)\ if $a_1=0$, $a_2\neq 0$, $a_4\neq 0$, then $\displaystyle f=\frac{a_4}{a_2}\ln(a_2\ln x+a_3)+C$,\\
where $C$ is an arbitrary constant.

Considering solutions 1)--4), we have the following possible forms of the function $f(x)$:
\begin{gather}\nonumber
f(x)=k_1x^n+k_2 , \ f=k_1(\ln x+k_2)^n+k_3  , \ \text{or} \ f=k_1\ln(\ln x+k_2)+k_3
\end{gather}
where $n$, $k_i$, $i\in\{1,2,3\}$ are arbitrary constants and $n\neq 0$, $k_1\neq 0$.

Using the set of equivalence transformations \eqref{Sp-St-13} we get a simplified view of
these functions in equation \eqref{Sp-St-3} (see {\bf Table 1}).

\centerline{{\it  {\bf Table~1.}} Classification of function $f(x)$}
\begin{center}
\begin{tabular}{|@{\;}l@{\;}|@{\;}c@{\;}|@{\;}c@{\;}|@{\;}c@{\;}|}  \hline
&&&\\[-3mm]
№ No& $f(x)$ &Equivalence transformations & $\overline{f}(\overline{x})$  \\[0.5mm] \hline
&&&\\[-3mm]
1& $k_1x^n+k_2$, &$\overline{t}=n^2t$, $\overline{x}=x^n$, $\displaystyle\overline{y}=\frac{n^2}{k_1}(k_2t+y)$,& $\overline{x}$\\
 &$k_1n\neq 0$& $\overline{u}=e^{(1-n^2)t/4}x^{(n-1)/2}u$&  \\[0.5mm] \hline
&&&\\[-3mm]
2  & $k_1(\ln x+k_2)^n+k_3$, &$\overline{t}=t$, $\overline{x}=e^{k_2}x$, & $\ln^n\overline{x}$ \\
 & $k_1n\neq 0$& $\displaystyle \overline{y}=\frac{1}{k_1}(k_3t+y)$,\ $\overline{u}=u$& \\[0.5mm] \hline
&&&\\[-3mm]
3 &$k_1\ln(\ln x+k_2)+k_3$, &$\overline{t}=t$, $\overline{x}=e^{k_2}x$, & $\ln\ln\overline{x}$\\
 &$k_1\neq 0$&$\displaystyle \overline{y}=\frac{1}{k_1}(k_3t+y)$,\ $\overline{u}=u$& \\[0.5mm] \hline
\end{tabular}
\end{center}

The equations with the function $f(x)$ from the list
\begin{gather}\label{Sp-St-17}
f(x)=x,\ \ f(x)=\ln^n x,\ \ f(x)=\ln\ln x,
\end{gather}
are mutually inequivalent with respect to the transformations from $G^{equiv}$. Thus, we obtain the following statement.

\begin{theorem}\label{thm:Number3}
If an equation of the form \eqref{Sp-St-3} admits algebra $A^{max}$, then the function $f(x)$ is equivalent with respect
to the transformations $G^{equiv}$ one of the forms $x$, $\ln^n x$, $\ln\ln x$, where $n\neq 0$.
\end{theorem}

The last step of algorithm of group classification is finding all possible Lie symmetries of each of those equations of class \eqref{Sp-St-3}
with the function $f(x)$ from the list \eqref{Sp-St-17}. For this, we substitute each function \eqref{Sp-St-17} to equations
\eqref{Sp-St-6}--\eqref{Sp-St-9} and get the general form of the components of a Lie symmetry operator of an equation from class \eqref{Sp-St-3}.

{\bf Case (1):\ $f(x)=x$}.

Substituting function $f(x)=x$ in equations \eqref{Sp-St-6}--\eqref{Sp-St-9}, we obtain
\begin{gather*}
\displaystyle \xi^0=C_1,\ \xi^1=(C_2y+C_3)x,\ \xi^2=\frac{1}{2}C_2y^2+C_3y+C_4, \\
\displaystyle \eta=\left(\frac{1}{2}C_2x+C_5\right)u+\beta(t,x,y),
\end{gather*}
where $C_i$, $i\in\{1,\ldots,5\}$ are arbitrary constants, $\beta(t,x,y)$ is an arbitrary solution of equation \eqref{Sp-St-3}.

{\bf Case (2):\ $f(x)=\ln^n x ,$ $n\neq -2, 0, 1$}.

Substituting function $f(x)=\ln^n x ,$ $n\neq -2, 0, 1$ in equations \eqref{Sp-St-6}--\eqref{Sp-St-9}, we obtain
\begin{gather*}
\displaystyle \xi^0=C_1t+C_2,\ \xi^1=\frac{1}{2}C_1x\ln x,\ \xi^2=\frac{n+2}{2}C_1y+C_3, \\
\displaystyle \eta=\left(\frac{1}{4}C_1(\ln x-t)+C_4\right)u+\beta(t,x,y),
\end{gather*}
where $C_i$, $i\in\{1,\ldots,4\}$ are arbitrary constants, $\beta(t,x,y)$ is an arbitrary solution of equation \eqref{Sp-St-3}.

{\bf Case (3):\ $f(x)=\ln x$}.

Substituting function $f(x)=\ln x$ in equations \eqref{Sp-St-6}--\eqref{Sp-St-9}, we obtain
\begin{gather*}
\xi^0=C_1t^2+C_2t+C_3,\\
\displaystyle \xi^1=\left(C_1t+\frac{1}{2}C_2\right)x\ln x+(C_4t^2+C_5t-3C_1y+C_6)x,\\
\displaystyle \xi^2=3C_1ty+\frac{3}{2}C_2y-\frac{1}{3}C_4t^3-\frac{1}{2}C_5t^2-C_6t+C_7, \\
\displaystyle \eta=\left(-C_1\ln^2 x+\frac{1}{2}\left((C_1-2C_4)t-C_5+\frac{1}{2}C_2\right)\ln x+\frac{1}{4}(2C_4-C_1)t^2+\right.\\
\displaystyle \left.+\frac{1}{4}(2C_5-8C_1-C_2)t-\left(\frac{3}{2}C_1+C_4\right)y+C_8\right)u+\beta(t,x,y) ,
\end{gather*}
where $C_i$, $i\in\{1,\ldots,8\}$ are arbitrary constants, $\beta(t,x,y)$ is an arbitrary solution of equation \eqref{Sp-St-3}.

{\bf Case (4):\ $f(x)=\ln^{-2} x$}.

Substituting function $f(x)=\ln^{-2} x$ in equations \eqref{Sp-St-6}--\eqref{Sp-St-9}, we obtain
\begin{gather*}
\xi^0=C_1t^2+C_2t+C_3,\ \ \displaystyle \xi^1=\left(C_1t+\frac{1}{2}C_2\right)x\ln x,\ \ \displaystyle \xi^2=C_4, \\
\displaystyle \eta=\left(-\frac{1}{4}C_1\ln^2 x+\frac{2C_1t+C_2}{4}\ln x-\frac{1}{4}C_1t^2-\frac{2C_1+C_2}{4}t+C_5\right)u+\beta(t,x,y) ,
\end{gather*}
where $C_i$, $i\in\{1,\ldots,5\}$ are arbitrary constants, $\beta(t,x,y)$ is an arbitrary solution of equation \eqref{Sp-St-3}.

{\bf Case (5):\ $f(x)=\ln\ln x$}.

Substituting function $f(x)=\ln\ln x$ in equations \eqref{Sp-St-6}--\eqref{Sp-St-9}, we obtain
\begin{gather*}
\displaystyle \xi^0=C_1t+C_2,\ \xi^1=\frac{1}{2}C_1x\ln x,\ \xi^2=C_1\left(y-\frac{1}{2}t\right)+C_3, \\
\displaystyle \eta=\left(\frac{1}{4}C_1(\ln x-t)+C_4\right)u+\beta(t,x,y),
\end{gather*}
where $C_i$, $i\in\{1,\ldots,4\}$ are arbitrary constants, $\beta(t,x,y)$ is an arbitrary solution of equation \eqref{Sp-St-3}.

The results are summarized in the following assertion.

\begin{theorem}\label{thm:Number4}
All possible equations of the form \eqref{Sp-St-3} admitting Lie algebras $A^{max}\supset A^{ker}$ of symmetries
are reduced to one of the $5$ “canonical” equations with functions given in {\bf Table 1} by an equivalence transformation
from $G^{equiv}$ and the finite-dimensional part of the maximal algebras of invariance of these “canonical” equations \eqref{Sp-St-3}
are presented in the third column of {\bf Table 2}.
\end{theorem}

\newpage
\centerline{{\it  {\bf Table~2.}} Results of group classification of class \eqref{Sp-St-3} with respect to $G^{equiv}$}

\begin{center}
\begin{tabular}{|@{\;}l@{\;}|@{\;}c@{\;}|@{\;}c@{\;}|}  \hline
&&\\[-3mm]
No& $f(x)$ &Basis of $A^{max}$ \\[0.5mm] \hline
&&\\[-3mm]
1& $\forall$ &$\partial_t$, $\partial_y$, $u\partial_u$\\[0.5mm] \hline
&&\\[-3mm]
2  & $x$ &$\partial_t$, $\partial_y$, $u\partial_u$, $x\partial_x+y\partial_y$, \\
 & & $\displaystyle xy\partial_x+\frac{1}{2}y^2\partial_y+\frac{1}{2}xu\partial_u$\\[2mm] \hline
&&\\[-3mm]
3  & $\ln^n x$ &$\partial_t$, $\partial_y$, $u\partial_u$, \\[1mm]
 &$n\neq -2, 0, 1$ & $\displaystyle t\partial_t+\frac{1}{2}x\ln x\partial_x+\frac{n+2}{2}y\partial_y+\frac{1}{4}(\ln x-t)u\partial_u$\\[2mm] \hline
&&\\[-3mm]
4 & $\ln x$ &$\partial_t$, $\partial_y$, $u\partial_u$, $x\partial_x-t\partial_y$, \\[1mm]
 & & $\displaystyle tx\partial_x-\frac{1}{2}t^2\partial_y+\frac{1}{2}(t-\ln x)u\partial_u$,\\[3mm]
 & & $\displaystyle t^2x\partial_x-\frac{1}{3}t^3\partial_y-\left(t\ln x+y-\frac{1}{2}t^2\right)u\partial_u$,\\[3mm]
& & $\displaystyle t\partial_t+\frac{1}{2}x\ln x\partial_x+\frac{3}{2}y\partial_y+\frac{1}{4}(\ln x-t)u\partial_u$, \\[3mm]
& & $\displaystyle t^2\partial_t+(t\ln x-3y)x\partial_x+3ty\partial_y-$ \\[3mm]
& & $\displaystyle -\left(\ln^2 x-\frac{1}{2}t\ln x+\frac{1}{4}t^2+2t+\frac{3}{2}y\right)u\partial_u$ \\[3mm] \hline
&&\\[-3mm]
5 & $\ln^{-2} x$ &$\partial_t$, $\partial_y$, $u\partial_u$, \\[1mm]
 & & $\displaystyle t\partial_t+\frac{1}{2}x\ln x\partial_x+\frac{1}{4}(\ln x-t)u\partial_u$,\\[3mm]
 & & $\displaystyle t^2\partial_t+tx\ln x\partial_x-$\\[3mm]
& & $\displaystyle -\frac{1}{4}(\ln^2 x-2t\ln x+t^2+2t)u\partial_u$ \\[3mm] \hline
&&\\[-3mm]
6 &$\ln\ln x$, &$\partial_t$, $\partial_y$, $u\partial_u$, \\
 & &$\displaystyle t\partial_t+\frac{1}{2}x\ln x\partial_x+\left(y-\frac{1}{2}t\right)\partial_y+ \frac{1}{4}(\ln x-t)u\partial_u$\\[0.5mm] \hline
\end{tabular}
\end{center}

In {\bf Table 2} we will not take into account the symmetry operator $X=\beta(t,x,y)\partial_u$ where $\beta(t,x,y)$ is an arbitrary solution of corresponding equation with function $f(x)$ which is inherent in linear equation and determines the principle of superposition.
The task of describing such operators is equivalent to the search for the general solution of such equations.

\end{document}